\documentclass{article}
\usepackage{authblk}
\usepackage{geometry}
\usepackage[titletoc]{appendix}
\usepackage{enumitem}
\usepackage{booktabs,diagbox}
\usepackage{hyperref}
\usepackage{xcolor}

\usepackage{tikz} 
\usetikzlibrary{cd}
\usetikzlibrary{arrows}

\usepackage{amssymb,amsmath,amsthm,mathtools,slashed,bm} 
\def\-{\raisebox{.75pt}{-}}

\numberwithin{table}{section}
\numberwithin{figure}{section}
\numberwithin{equation}{section}

\theoremstyle{definition}
\newtheorem{defn}{Definition}[section]

\newtheorem{exmp}{Example}[section]
\newtheorem{rmk}{Remark}[section]
\theoremstyle{plain}
\newtheorem{lem}{Lemma}[section]
\newtheorem{prop}{Proposition}[section]
\newtheorem{thm}{Theorem}[section]
\newtheorem*{thm*}{Theorem}

\usepackage{amsrefs}

\title{\textbf{Gauge Theory on Graphs}}
\author[1]{Shuhan Jiang}
\date{}
\affil[1]{Max Planck Institute for Mathematics in the Sciences, 04103 Leipzig}

\begin{document}
	
	\maketitle
	
	\begin{abstract}
		In this paper, we provide the notions of connection $1$-forms and curvature $2$-forms on graphs. We prove a Weitzenböck formula for connection Laplacians in this setting. We also define a discrete Yang-Mills functional and study its Euler-Lagrange equations.
	\end{abstract}
	%\tableofcontents
	\section{Introduction}
	
	Gauge theory is of fundamental importance in both modern physics and mathematics. From a physical perspective, it provides an elegant geometric framework to describe electromagnetism, weak interactions, and strong interactions uniformly. From a mathematical perspective, the study of the Euler-Lagrange equations of gauge invariant functionals reveals astounding information about smooth structures of low dimensional manifolds \cites{Donaldson1997,Morgan1996}. 
	
	To understand the quantum aspects of gauge theory better, the physicist Kenneth Wilson formulated it on lattices in 1974 \cite{Wilson1974}. His methods were later developed into a mature area of research nowadays known as lattice gauge theory \cite{Wilson2005}. Lattice gauge theory has many successful applications in quantum chromodynamics calculations \cite{Gattringer2009}. However, the geometry flavor of gauge theory is lost after this discretization. On the other hand, there exist other discretizations of physical theories where the geometry flavors of the original theories are preserved, for example, Robin Forman's discretization \cite{Forman1998} of Witten-Morse theory \cite{Witten1982}. Inspired by Forman's work, here we present a formulation of gauge theory on a general graph. 
	
	The paper is organized as follows:
	
	 In Section 2, we review some basic concepts in differential geometry and graph theory and provide a slight generalization of the Weitzenböck decomposition of a symmetric matrix \cite{Forman2003}. In Section 3, we present the standard definitions of differential forms, exterior derivatives, and (Hodge) Laplacians on graphs. Additionally, we briefly recall the Weitzenböck formula for Laplacians \cite{Forman1998}.
	 In Section 4, we extend the notion of an exterior derivative $d$ to an exterior covariant derivative $d_A$, introduce the curvature $2$-form $F$ defined through the standard formula $F=d_A \circ d_A$, and demonstrate that $F$ satisfies the second Bianchi identity $d_A F=0$.
	 In Section 5, we extend the notion of a Laplacian $\Delta$ to a connection Laplacian $\Delta_A$ and establish a Weitzenböck formula in this setting, which generalizes Forman's result in \cite{Forman1998}.
	 In Section 6, we define a Yang-Mills functional and derive its Euler-Lagrange equations. We demonstrate that the minimum and maximum of the functional can always be achieved and provide several examples for the abelian gauge group $\mathrm{U}(1)$.
	 Finally, in Section 7, we discuss potential future research directions stemming from the present work.
	
	The research for this paper was performed in 2020. Later, the paper by \cite{Bianconi2021} by Ginestra Bianconi appeared where she also presented an approach to gauge theory on graphs. More recently, she has further developed her approach. It turns out, however, that her and my approach are different from each other.
	
	\section{Preliminaries}
	
	In this section, we review some basic definitions in differential geometry and graph theory. For more details, we refer the reader to \cite{Kobayashi1963,Jost2017,Rudolph2017,Hatcher2002}.
	
	\subsection{Geometry preliminaries}
	
	\begin{defn}
		A vector bundle of rank $n$ is a triple $(E,M,\pi)$ where $E$ and $M$ are smooth manifolds, $\pi: E \rightarrow M$ is a surjective smooth map, such that for any $x \in M$, the fiber $E_x:=\pi^{-1}(x)$ carries the structure of an $n$-dimensional vector space, and there exists an open neighborhood $U$ of $x$ and a diffeomorphism $\varphi: \pi^{-1}(U) \rightarrow U \times \mathbb{R}^n$ of the form $\varphi(p)=(\pi(p),\psi(p))$, where $\psi|_{\pi^{-1}(y)}$ is a vector space isomorphism for any $y \in U$. 
	\end{defn}
    The product space $M \times \mathbb{R}^n$ can be viewed as a vector bundle, which we call the trivial vector bundle of rank $n$ over $M$. For any two vector bundle $E$ and $F$ over a common base manifold $M$. The disjoint unions
    $E \oplus F:=\coprod_{x \in M} E_x \oplus F_x$ and $E \otimes F:=\coprod_{x \in M} E_x \otimes_{\mathbb{R}} F_x$ can be given vector bundle structures and are called the sum and the tensor product of $E$ and $F$, respectively. A homomorphism between $E$ and $F$ is a smooth map $f: E \rightarrow F$ such that $\pi_E = \pi_F \circ f$ and $f$ is a fiber-wise linear map. The disjoint union
    $\mathrm{Hom}(E,F):=\coprod_{x \in M} \mathrm{Hom}_{\mathbb{R}}(E_x, F_x)$ can be also given a vector bundle structure. In particular, we call $\mathrm{End}(E):=\mathrm{Hom}(E,E)$ the endomorphism bundle of $E$. Another important example is the exterior bundle $\Lambda E$ of $E$, whose fiber at $x \in M$ is the exterior algebra $\Lambda(E_x)$.
    % , and $E^*= \mathrm{Hom}(E,M \times \mathbb{R})$ the dual bundle of $E$.
    \begin{defn}
    	A section $s$ of $E$ is a smooth map $s: M \rightarrow E$ such that $\pi \circ s = \mathrm{id}_M$. The space of sections of $E$ is denoted by $\Gamma(E)$.
    \end{defn}
    $\Gamma(E)$ is an infinite dimensional vector space and is canonically a $C^{\infty}(M)$-module, where $C^{\infty}(M)$ is the algebra of smooth functions over $M$. Note that $\Gamma(E \otimes F)=\Gamma(E) \otimes_{\mathbb{R}} \Gamma(F)$.
    \begin{defn}
    	A covariant derivative $\nabla$ is a $\mathbb{R}$-linear map $\nabla: \Gamma(E) \rightarrow \Gamma(T^*M \otimes E)$ such that
    	\begin{align*}
    		\nabla(fs)=d(f) \otimes_{\mathbb{R}} s + f \otimes_{\mathbb{R}} \nabla(s),
    	\end{align*}
        for any $f \in C^{\infty}(M)$ and $s \in \Gamma(E)$, where $d$ is the de Rham differential of $M$. If $E$ carries a metric, i.e., if $E_x$ has an inner product for all $x \in M$, we have an induced symmetric bilinear pairing $\langle \cdot,\cdot \rangle: \Gamma(E) \times \Gamma(E) \rightarrow C^{\infty}(M)$. $\nabla$ is said to be compatible with the metric of $E$ if 
        \begin{align*}
        	d\langle s,t \rangle = \langle \nabla(s),t \rangle + \langle s,\nabla(t) \rangle
        \end{align*}
        for all $s, t \in \Gamma(E)$.
    \end{defn}
    Note that $\Gamma(E)$ is also canonically a $\Gamma(\mathrm{End}(E))$-module. In fact, the difference between two covariant derivatives can be seen as a section of $\Omega^1(M) \otimes \Gamma(\mathrm{End}(E))$. A covariant derivative $\nabla$ of $E$ induces a covariant derivative of the endomorphism bundle $\mathrm{End}(E)$, which we again denote by $\nabla$, by requiring
    \begin{align*}
    	\nabla(H(s)) = \nabla(H)(s) + H(\nabla(s))
    \end{align*}
    for any $H \in \Gamma(\mathrm{End}(E))$ and $s \in \Gamma(E)$. 
    
    Let $\Omega^p(M,E):= \Gamma(\Lambda^p T^*M \otimes E)$. $\nabla$ can be easily extended to a map $d_{\nabla}: \Omega^p(M,E) \rightarrow \Omega^{p+1}(M,E)$ satisfying 
    \begin{align*}
    	d_{\nabla} (\omega \otimes_{\mathbb{R}} s) = d \omega \otimes_{\mathbb{R}} s + \omega \wedge \nabla(s),
    \end{align*}
    where $\omega$ is a $p$-form over $M$ and $\wedge$ is the wedge product of differential forms. 
    \begin{prop}
    	There exists a $2$-form $F$ with values in $\mathrm{End}(E)$ such that
    	\begin{align*}
    		F = d_{\nabla} \circ d_{\nabla}.
    	\end{align*} 
        $F$ is called the curvature of $\nabla$. $F$ satisfies the so-called second Bianchi identity
        \begin{align*}
        	d_{\nabla} F = 0.
        \end{align*}
    \end{prop}

	Let $G$ be a Lie group. 
	\begin{defn}
		A principal $G$-bundle is a triple $(P,M,\pi)$ where $P$ and $M$ are smooth manifolds, $\pi: P \rightarrow M$ is a surjective smooth map, with a free proper smooth $G$-action on $P$, such that $\pi(pg)=\pi(p)$ for any $p \in P$ and $g \in G$, and for every $x \in M$, there exists an open neighborhood $U$ of $x$ and a diffeomorphism $\varphi: \pi^{-1}(U) \rightarrow U \times G$ of the form $\varphi(p)=(\pi(p),\psi(p))$ which is $G$-equivariant, i.e., $\varphi(pg)=(\pi(p),\psi(p)g)$ for all $g \in G$.
	\end{defn}
    Let $\mathfrak{g}$ be the Lie algebra of $G$. A differential form $\omega$ over $P$ is called horizontal if $\iota_{X_{\xi}} \omega = 0$ for all $\xi \in \mathfrak{g}$, where $X_{\xi}$ is the fundamental vector field of $\xi$, and $\iota_{X_{\xi}}$ is the contraction by $X_{\xi}$. A horizontal differential form $\omega$ is called basic if $\mathrm{Lie}_{X_{\xi}} \omega = 0$ for all $\xi \in \mathfrak{g}$, where $\mathrm{Lie}_{X_{\xi}}$ is the Lie derivative induced by $X_{\xi}$. 
    \begin{lem}\label{bijform}
    	There exists a bijection between the set of basic differential forms over $P$ and the set of differential forms over $M$.
    \end{lem}
    
    A $\mathfrak{g}$-valued differential form over $P$ is an element in $\Omega(P) \otimes_{\mathbb{R}} \mathfrak{g}$.
    
    \begin{defn}
    	A $\mathfrak{g}$-valued $1$-form $A$ over $P$ is called a connection ($1$-form) if 
    	\begin{align*}
    		A(X_{\xi}) = \xi
    	\end{align*}
        for all $\xi \in \mathfrak{g}$. The curvature ($2$-form) of $A$ is a  $\mathfrak{g}$-valued $2$-form $F:= dA + \frac{1}{2}[A,A]$.
    \end{defn}
            
    Let $\rho: G \rightarrow \mathrm{GL}(V)$ be a representation of $G$ on a vector space $V$. Define the $G$-action on $P \times V$ via $(p,v)g=(pg,\rho(g^{-1})v)$. Let $P\times_{\rho}V$ denote the quotient $(P \times V)/G$. There exists a well-defined projection $P \times_{\rho} V \rightarrow M$ sending $[p,v]$ to $\pi(p)$, which we will still denote by $\pi$. The triple $(P \times_{\rho} V,M,\pi)$ is a vector bundle over $M$.
    \begin{defn}
    	$(P \times_{\rho} V,M,\pi)$ is called the associated bundle of $(P,M,\pi)$.
    \end{defn}
    %Conversely, given a vector bundle $E$, the disjoint union of frames  
    In particular, let $\rho$ be the adjoint representation of $G$ on $\mathfrak{g}$. The associated bundle corresponding to $\rho$ is called the adjoint bundle of $P$, denoted by $\mathrm{ad} P$. Let $\mathcal{A}$ denote the space of connections on $P$. Note that the difference between any two connections is a $G$-invariant basic form on $P$. Therefore, $\mathcal{A}$ is an affine space modelled on the space of sections of the vector bundle $\Omega^1(M) \otimes \mathrm{ad} P$.

    \begin{prop}\label{spconn}
    	Let $P$ be a principal $\mathrm{GL}(\mathbb{R}^n)$-bundle. Let $E$ be a vector bundle of rank $n$ associated to $P$ through the fundamental representation of $\mathrm{GL}(\mathbb{R}^n)$. There exists a bijection between the set of connections on $P$ and the set of covariant derivatives of $E$. 
    	
    	Moreover, let $\nabla_A$ denote the covariant derivative induced by $A$. The curvature of $A$ is a basic $2$-form on $P$ and can be identified with the curvature of $\nabla_A$ using Lemma \ref{bijform}.
    \end{prop}
    %Given a vector bundle $E$,
    \begin{rmk}
    	With a slight abuse of notation, we will denote the exterior covariant derivative associated to $A$ by $d_A$ instead of $d_{\nabla_A}$. We will not distinguish between the words ''(exterior) covariant derivative" and ''connection". We will simply call a connection which is compatible with the metric a metric connection.
    \end{rmk} 
    \begin{rmk}\label{spconnad}
    	Proposition \ref{spconn} can be refined as follows. If the fiber of $E$ carries an additional structure, for example, an inner product, then there exists a subbundle $P_o$ of $P$, which is a principal $\mathrm{O}(n)$-bundle, and a bijection between the set of connections on $P_o$ and the set of metric connections of $E$. If $n=2m$ and the fiber of $E$ carries a complex structure, i.e., $E_x \cong \mathbb{C}^m$, $x \in M$, then there exists a subbundle $P_c$ of $P$, which is a principal $\mathrm{GL}(\mathbb{C}^m)$-bundle, and a bijection between the set of connections on $P_c$ and the set of complex linear connections of $E$.  If $E_x \cong \mathbb{C}^m$ and carries a Hermitian inner product, $x \in M$, then there exists a subbundle $P_u$ of $P$, which is a principal $\mathrm{U}(m)$-bundle, and a bijection between the set of connections on $P_u$ and the set of complex linear metric connections of $E$.
    	%Other structures such as orientations and Hermitian inner products can be also considered.
    \end{rmk}
    
    \begin{defn}
    	The automorphism group
    	\begin{align*}
    		\mathcal{G}:=\{f: P \rightarrow P | \pi \circ f = \pi, ~ f(pg)=f(p)g, ~ \forall p \in P, g \in G \}
    	\end{align*}
        of $P$ is called the gauge symmetry group. Elements of $\mathcal{G}$ are called gauge transformations.
    \end{defn}
	$\mathcal{G}$ acts naturally on $\mathcal{A}$ via pullbacks. If $P$ is trivial, $\mathcal{G}$ is nothing but the mapping space $C^{\infty}(M,G)$ with the canonical group structure, and $\mathcal{A} \cong \Omega^1(M) \otimes_{\mathbb{R}} \mathfrak{g}$. 
	%Let $E$ be a vector bundle associated to $P$ via $\rho: G \rightarrow \mathrm{Gl}(V)$.
	$g \in \mathcal{G}$ and $A \in \mathcal{A}$ can be viewed as $\mathrm{GL}(\mathbb{R}^n)$-valued and $\mathrm{End}(\mathbb{R}^n)$-valued functions over $M$, respectively. The action of $\mathcal{G}$ on $\mathcal{A}$ is then given by
	\begin{align*}
		%gA = gdg^{-1} + gAg^{-1}.
		(gA)(x) = g(x)(dg^{-1})(x) + g(x)A(x)g^{-1}(x).
	\end{align*}
    %Locally, 
    $\nabla_A$ can be written as $d+A$. It is easy to check that
    \begin{align*}
    	\nabla_{g A} = g \circ \nabla_{A} \circ g^{-1}.
    \end{align*}
    It follows that
    \begin{align*}
    	gF = g \circ F \circ g^{-1}.
    \end{align*}
	Now let $\gamma: [0,1] \rightarrow M$ be a smooth curve in $M$ with $\gamma(0)=x$ and $\gamma(1)=y$. Let $s \in \Gamma(E)$. We say $s$ is parallel along $\gamma$ with respect to $A$ if it is a solution of
	$(\nabla_A s (t))(\dot{\gamma}(t))=0$. Such a solution is uniquely determined by the initial value $s(0)$ by the main theorem of ODE. Therefore, $s$ determines a bijective map from $E_x$ to $E_y$ which we denote by $A(y,x)$. A gauge transformation $g$ acts on $A(y,x)$ as
	\begin{align*}
		(gA)(y,x)=g(y)A(y,x)g(x)^{-1}.
	\end{align*}
	
	Let $G=\mathrm{O}(n)$ (or $\mathrm{U}(n)$). Let $E$ be a vector bundle associated to $P$ via the standard representation $\rho: G \rightarrow \mathrm{GL}(\mathbb{R}^n)$ (or $\rho: G \rightarrow \mathrm{GL}(\mathbb{C}^n$)). The standard inner product on $\mathbb{R}^n$ (or $\mathbb{C}^n$) induces a $G$-invariant metric $\langle \cdot,\cdot \rangle$ on $E$. By Proposition \ref{spconn} and Remark \ref{spconnad}, $\mathcal{A}$ can also be viewed as the space of (complex linear) metric connections of $E$. Let $M$ be equipped with a Riemannian metric $g$. $g$ and $\langle \cdot,\cdot \rangle$ induce a metric on the vector bundle $\Lambda T^* M \otimes E$, which we denote again by $\langle \cdot,\cdot \rangle$. Let $\mathrm{dvol}$ be the volume form on $M$.
	
	\begin{defn}
		The Yang-Mills functional is the $\mathcal{G}$-invariant functional on $\mathcal{A}$ defined by
		\begin{align*}
			\mathcal{YM}(A) = \int_M \mathrm{dvol} \langle F,F \rangle.
		\end{align*}
	\end{defn}
	Let $d_A^*$ denote the adjoint of $d_A$ with respect to $\int_M \mathrm{dvol} \langle \cdot,\cdot \rangle$. The Euler-Lagrange equations of the Yang-Mills functional are
	\begin{align*}
		d_A^* F = 0. 
	\end{align*}
    %\begin{rmk}
    %	\textcolor{red}{Motivation.}
    %\end{rmk}
	From now on, we use $\nabla$ to denote the Levi-Civita connection on $\Lambda T^*M$ and use $\nabla^*$ to denote its adjoint with respect to the Riemannian metric $g$ on $M$. Recall that for a differential form $T$, we have
	\begin{align}\label{wb}
		\Delta T = \nabla^*\nabla T + \mathrm{Ric}(T),
	\end{align}
	where $\Delta$ is the Hodge Laplacian and $\mathrm{Ric}(\cdot)$ is the so-called Weitzenböck operator defined by
	\begin{align}\label{wo}
		\mathrm{Ric}(T)(X_1,\cdots,X_k) = \sum_{i=1}^n \sum_{j=1}^k (R(e_i,X_j)T)(X_1,\cdots,e_i,\cdots,X_k),
	\end{align}
	where $R$ is the Riemannian curvature tensor and $\{e_i\}_{i=1}^n$ is a local orthonormal frame. Applying \eqref{wo} to the metric dual $X^{\sharp}=g(X,\cdot)$ of the vector field $X$, we have
	\begin{align*}
		\mathrm{Ric}(X^{\sharp})(Y)=\sum_{i=1}^n (R(e_i,Y)X^{\sharp})(e_i)
		=\sum_{i=1}^n g(R(e_i,Y)X,e_i)
		=-\sum_{i=1}^n g(R(e_i,Y)e_i,X)
		=\mathrm{Ric}(X,Y),
	\end{align*} 
	where $\mathrm{Ric}(\cdot,\cdot)$ is the Ricci curvature tensor of $g$, hence the notation. 
	
	Let $E$ be a vector bundle over $M$ equipped with a metric. Let $\nabla_A$ be a metric connection of $E$, $\nabla_A$ and the Levi-Civita connection induces a connection on the tensor product $\Lambda T^*M \otimes E$, which we denote by $\nabla_A$ with a slight abuse of notation. \eqref{wb} can be generalized to the case of vector valued differential forms. Let $F$ be the curvature $2$-form of $A$. For $T \in \Gamma(\Lambda T^*M \otimes E)$, define similarly
	\begin{align*}
		F(T)(X_1,\cdots,X_k) = \sum_{i=1}^n \sum_{j=1}^k (F(e_i,X_j)T)(X_1,\cdots,e_i,\cdots,X_k).
	\end{align*}
	The generalized Weitzenböck formula takes the form \cite{Bourguignon1981}
	\begin{align}\label{wbg}
		\Delta_A T = \nabla_A^*\nabla_A T + \mathrm{Ric}(T) + F(T)
	\end{align}
	where $\Delta_A = d_Ad_A^*+d_A^*d_A$, and $\nabla_A^*$ is the adjoint of $\nabla_A$ with respect to the inner product on $\Gamma(\Lambda T^*M \otimes E)$ induced by the Riemannian metric of $M$ and the metric of $E$.
	
	\subsection{Graph theory preliminaries}
	
	\begin{defn}
		A (finite simple) graph is a pair $\Gamma=(V,E)$, where $V$ is a finite set, $E \subseteq \binom{V}{2}$ is a set of $2$-subsets of $V$. The elements of $V$ are called vertices, and the elements of $E$ are called edges. 
	\end{defn}

	Once the set $E$ is specified, we automatically get sets of ''cliques" of higher order. 
	\begin{defn}
		Let $k$ be a natural number. A $k$-clique is a $k$-subset $\sigma$ of $V$ such that for all distinct $i, j \in \sigma$, $\{i,j\} \in E$.  
	\end{defn}
	In other words, a $k$-clique is just a complete $k$-subgraph of $\Gamma$. We use $K_k(\Gamma)$ to denote the set of $k$-cliques of $\Gamma$. A maximum clique is a clique of maximum possible size in $\Gamma$. Its size is known as the clique number of $\Gamma$, denoted by $\omega(\Gamma)$. 
	%(From now on, we will always assume that $\omega(\Gamma)<\infty$.) 
	Let $K(\Gamma)$ be the set of all cliques of a graph $\Gamma$, i.e.,
	\[
	K(\Gamma) = \coprod_{k=1}^{\omega(\Gamma)} K_k(\Gamma).
	\]
	
	\begin{defn}
		An abstract simplicial complex is a collection $X$ of non-empty finite subsets of a set $S$ such that for every $\alpha \in X$ and every non-empty subset $\beta \subset \alpha$, $\beta \in X$. 
	\end{defn}

    The vertex set $V(X)$ of $X$ is defined to be the union $\bigcup_{\alpha \in X} \alpha$. Elements of $X$ are called simplices. For two simplices $\alpha$ and $\beta$, we write $\alpha < \beta$ or $\beta > \alpha$ if $\alpha \subset \beta$. 
    
    As an example, $K(\Gamma)$ is an abstract simplicial complex and is called as the clique complex of $\Gamma$. Note that a $(k+1)$-clique is a $k$-simplex.
    
    \begin{rmk}[Caveat]
    	$\Gamma$ is typically viewed as a $1$-dimensional CW complex, i.e., the geometric realization of the abstract simplicial complex $\coprod_{k=1}^2 K_k(\Gamma)$. Our convention of the topological space associated to $\Gamma$ is, however, a $(\omega(\Gamma)-1)$-dimensional CW complex, i.e., the geometric realization of its clique complex $K(\Gamma)$.
    \end{rmk}
	
	\begin{defn}
		A graph homomorphism between two graphs $m: \Gamma_1=(V_1,E_1) \rightarrow \Gamma_2=(V_2,E_2)$ is a map from $V_1$ to $V_2$ such that $\{m(u),m(v)\} \in E_2$ if $\{u,v\} \in E_1$.		
	\end{defn}

    \begin{defn}
    	A simplicial map between two abstract simplicial complexes $m: X_1 \rightarrow X_2$ is a map from $V(X_1)$ to $V(X_2)$ such that $m(\alpha) \in X_2$ if $\alpha \in X_1$.		
    \end{defn}

    A graph homomorphism $m: \Gamma_1 \rightarrow \Gamma_2$ induces naturally a simplicial map $K(m)$ from $K(\Gamma_1)$ to $K(\Gamma_2)$. Let $m: \Gamma_1 \rightarrow \Gamma_2$ and $m': \Gamma_2 \rightarrow \Gamma_3$ be two graph homomorphisms. It is also easy to check that $K(m'\circ m)=K(m')\circ K(m)$. In other words, we obtain a functor $K$ from the category of (simple) graphs to the category of abstract simplicial complexes.
    
	\begin{defn}
		Let $\Gamma=(V,E)$ be a graph. An orientation on $\Gamma$ is an assignment of a direction to each edge $e \in E$, making $e$ into an ordered pair. Given an orientation,  a vertex sequence $\{v_1, v_2, \dots, v_{n+1}\}$ with $\{v_i, v_{i+1}\} \in E$, $v_i \rightarrow v_{i+1}$ for $i=1,\cdots,n$ is called a directed path in $\Gamma$ if $v_1 \neq v_{n+1}$, and a directed cycle in $\Gamma$ if $v_1=v_{n+1}$. An orientation on $\Gamma$ is said to be acyclic if there is no directed cycle in $\Gamma$ with respect to it. A graph equipped with an acyclic orientation is called a directed acyclic graph.
	\end{defn}

    %Every graph has an acyclic orientation. 
    An acyclic orientation of $\Gamma=(V,E)$ induces a partial order $<$ on $V$ as follows. For any $i,j \in V$, $i<j$ if and only if there exists a directed path in $\Gamma$ from $i$ to $j$. Since there is no directed cycle in $\Gamma$, $<$ is indeed well-defined.
	
	\begin{defn}
		An abstract simplicial complex $X$ is said to be ordered if there is a partial order on $V(X)$, restricting to a total order on each simplex.
	\end{defn}
    
    Similarly, we have a canonical functor from the category of (simple) directed acyclic graphs to the category of ordered abstract simplicial complexes.
    
    \subsection{Weitzenböck decomposition}
    
    Let $(R,\|\cdot\|,1_R)$ be a unital Banach algebra. Let $A=\left(A_{i j}\right)$ be a symmetric $n \times n$ matrix over $R$.
    
    \begin{defn}
    	Let $B(A)=\left(B_{i j}\right)$ denote the symmetric $n \times n$ matrix whose entries are given by
    	
    	$$
    	B_{i j}= \begin{cases}A_{i j} & \text { if } i \neq j \\ \sum_{j \neq i}\|A_{i j}\| & \text { if } i=j\end{cases}
    	$$
    	where we identity $r \in \mathbb{R}_{\geq 0}$ with $r 1_R \in R$. Let $\mathrm{Ric}(A)=\left(\mathrm{Ric}_{i j}\right)$ denote the $n \times n$ diagonal matrix whose entries are given by    	
    	$$
    	\mathrm{Ric}_{i j}= \begin{cases}0 & \text { if } \quad i \neq j \\ A_{i j}-\sum_{j \neq i}\|A_{i j}\| & \text { if } \quad i=j\end{cases}
    	$$    	
    	We refer to $B(A)$ as the Bochner matrix associated to $A$ and $\mathrm{Ric}(A)$ as the curvature matrix associated to $A$. Moreover, we refer to the identity    	
    	$$
    	A=B(A)+\mathrm{Ric}(A)
    	$$   	
    	as the Weitzenböck decomposition of $A$.
    \end{defn} 
	
	\section{Differential forms and Laplacians}
	
	Let $\Gamma=(V,E)$ be a graph. %From now on, we will always assume $\omega(\Gamma)<\infty$.   
    \begin{defn}
    A $k$-form on $\Gamma$ is a map $f: V \underbrace{\times \cdots \times}_{k+1} V \rightarrow \mathbb{R}$ such that
    \begin{align*}
    f(i_{\sigma(0)},\cdots,i_{\sigma(k)}) = (-1)^{|\sigma|}f(i_0,\cdots,i_k)
    \end{align*}
    for $\lbrace i_0, \cdots, i_k \rbrace \in K_{k+1}(\Gamma)$, where $\sigma \in S_{k+1}$ and $|\sigma|$ is the sign of $\sigma$, and
    \begin{align*}
    f(i_0,\cdots,i_k) = 0
    \end{align*}
    for $\lbrace i_0, \cdots, i_k \rbrace \notin K_{k+1}(\Gamma)$.     
    \end{defn}

    We denote the space of $k$-forms by $\Omega^k(\Gamma)$, and put
    \[
    \Omega(\Gamma)=\bigoplus_{k=0}^{\omega(\Gamma)-1} \Omega^k(\Gamma).
    \]
    
    \begin{defn}
    The $k$-th exterior derivative $d_k: \Omega^k(\Gamma) \rightarrow \Omega^{k+1}(\Gamma)$ is defined by 
    \begin{align*}
    (d_k f)(i_0,\cdots,i_{k+1})=\sum_{j=0}^{k+1}(-1)^j f(i_0,\cdots,\hat{i_j},\cdots,i_{k+1}).
    \end{align*}
    \end{defn}
       
    By definition, $d_{k+1} \circ d_k =0$. The $k$-th cohomology group of $\Gamma$ can be then defined as the quotient
    $
    	H^k(\Gamma)=\mathrm{ker}(d_k)/\mathrm{im}(d_{k-1}).
    $
    We use $d$ to denote the exterior derivative on $\Omega(\Gamma)$.
    
    One can put an inner product $\langle \cdot, \cdot \rangle$ on $\Omega(\Gamma)$ by setting
	\begin{align*}
	\langle f,g \rangle = \sum_{\{i_0, \cdots, i_k\} \in K_{k+1}(\Gamma)} f(i_0,\cdots,i_k)g(i_0,\cdots,i_k),
	\end{align*}
	where $f$ and $g$ are two $k$-forms. $\langle \cdot, \cdot \rangle$ is well defined because 
	$
	f(i_{\sigma(0)},\cdots,i_{\sigma(k)})g(i_{\sigma(0)},\cdots,i_{\sigma(k)})=f(i_0,\cdots,i_k)g(i_0,\cdots,i_k).
	$
	
	$\Omega(\Gamma)$ equipped with $\langle \cdot, \cdot \rangle$ becomes a finite dimensional Hilbert space. We denote the adjoint of the exterior derivative $d$ by $d^*$. It is not hard to work out an explicit formula for $d^*$. For a $k$-form $f$, we have
	\begin{align*}
	d^*f(i_0,\cdots,i_{k-1})= \sum_{\{l,i_0,\cdots,i_{k-1}\} \in K_{k+1}(\Gamma)} f(l, i_0,\cdots,i_{k-1}). 
	\end{align*}
	Since $f$ is only nonzero if $\{l,i_0,\cdots,i_{k-1}\}$ is a $(k+1)$-clique, we can also simply write
	\begin{align*}
		d^*f(i_0,\cdots,i_{k-1})= \sum_{l} f(l, i_0,\cdots,i_{k-1}). 
	\end{align*}
	\begin{defn}
	The (Hodge) Laplacian on $\Omega(\Gamma)$ is defined as
	\begin{align*}
	\Delta = dd^* + d^*d.
	\end{align*}
	A $k$-form $\varphi \in \Omega^k(\Gamma)$ is said to be harmonic if 
	$
	\Delta \varphi = 0.
	$
	\end{defn}
	Since $\Delta$ preserve the degree of a differential form, the space of harmonic forms $\ker \Delta$ is canonically a graded vector space.
	\begin{thm}
		$H(\Gamma) \cong \mathrm{ker}(\Delta)$ as graded vector spaces.
	\end{thm}
	\begin{proof}
		The proof can be found in, for example, \cites{Lim2020,Horak2013}.
	\end{proof}

	One can work out an explicit expression for $\Delta$. For $k=0$, we have
	\begin{align}
		\Delta f(i)= d^*df(i)=\sum_{\{l,i\} \in E} df(l,i)=\sum_{\{l,i\} \in E} (f(i)-f(l))=\deg(i) f(i) -  \sum_{\{l,i\} \in E} f(l),
	\end{align}
	where $\deg(i) = \# \{l \in V | \{l,i\} \in E \}$ is the degree of the vertex $i \in V$.
	
	For $k>0$, we have
	\begin{align*}
	&d^*df(i_0,\cdots,i_k)=\sum_{\{l,i_0,\cdots,i_{k}\} \in K_{k+2}(\Gamma)} df(l,i_0,\cdots,i_k) \\
	&=\sum_{\{l,i_0,\cdots,i_{k}\} \in K_{k+2}(\Gamma)} \left(f(i_0,\cdots,i_k)-\sum_{j=0}^k(-1)^j f(l,i_0,\cdots,\widehat{i_j},\cdots,i_k) \right) \\
	&=\deg(i_0,\cdots,i_k)f(i_0,\cdots,i_k) - \sum_{\{l,i_0,\cdots,i_{k}\} \in K_{k+2}(\Gamma)}\sum_{j=0}^k(-1)^j f(l,i_0,\cdots,\widehat{i_j},\cdots,i_k),
	\end{align*}
	where $\mathrm{deg}(i_0,\cdots,i_k)$ is the number of $(k+2)$-cliques containing $i_0,\cdots,i_k$.  
	
	On the other hand, we have
	\begin{align*}
	&dd^*f(i_0,\cdots,i_k)=d^*f(i_1,\cdots,i_k)+\sum_{j=1}^k(-1)^j d^*f(i_0,\cdots,\widehat{i_j},\cdots,i_k) \\
	&=\sum_{\{l,i_1,\cdots,i_{k}\} \in K_{k+1}(\Gamma)} f(l,i_1,\cdots,i_k) + \sum_{j=1}^k(-1)^j \sum_{\{l,i_0,\cdots,\widehat{i_j},\cdots,i_{k}\} \in K_{k+1}(\Gamma)} f(l,i_0,\cdots,\widehat{i_j},\cdots,i_k) \\
	&=(k+1)f(i_0,\cdots,i_k) + \sum_{j=0}^k (-1)^j \sum_{\substack{l \neq i_j \\ \{l,i_0,\cdots,\widehat{i_j},\cdots,i_{k}\} \in K_{k+1}(\Gamma)}} f(l,i_0,\cdots,\widehat{i_j},\cdots,i_k).
	\end{align*}
	It follows that
	\begin{align}
	\begin{split}
	\Delta f(i_0,\cdots,i_k)&=\left(\mathrm{deg}(i_0,\cdots,i_k) + (k+1) \right) f(i_0,\cdots,i_k) \\
	&+\sum_{j=0}^k(-1)^j \sum_{\substack{\{l,i_0,\cdots,\widehat{i_j},\cdots,i_{k}\} \in K_{k+1}(\Gamma)\\ \{l,i_0,\cdots,i_{k}\} \notin K_{k+2}(\Gamma) \\ l \neq i_j}} f(l,i_0,\cdots,\widehat{i_j},\cdots,i_k).
	\end{split}
	\end{align}
	
	In \cite{Forman2003}, Forman proposed a discrete analogue of \eqref{wb}. Let $X$ be a CW complex. If $\alpha$ and $\beta$ are cells of $X$, we write $\alpha<\beta$ or $\beta>\alpha$ if $\alpha$ is contained in the boundary of $\beta$. Two $k$-cells $\alpha$ and $\beta$ are called parallel neighbors if they share either a $(k+1)$-cell or $(k-1)$-cell but not both. The Forman-Ricci curvature of $X$ is a function $\mathcal{R}: X \rightarrow \mathbb{R}$ defined by
		\[
		\mathcal{R}(\alpha)=\#\{(k+1)\textrm{-cells}~\beta> \alpha \} + \#\{(k-1)\textrm{-cells}~\alpha < \beta \} - \#\{\textrm{parallel neighbors}\}.
		\]
	Note that for a simplicial complex, $\#\{(k-1)\textrm{-cells}~\alpha < \beta \} = k+1$.
	
	Using Forman's convention, we can write down a simpler expression for the Laplacian of $\Gamma$. For a $k$-form $f$ and a $k$-simplex (or a $(k+1)$-clique) $\alpha$, we have
	\begin{align}
		%&d f(\alpha) = \sum_{\gamma<\alpha} \epsilon_{\gamma\alpha} f(\gamma), \quad d^*f(\alpha) = \sum_{\gamma>\alpha} \epsilon_{\alpha\gamma} f(\gamma) \\
		&\Delta f(\alpha)=\left(\#\{(k+1)\textrm{-simplices}~\beta> \alpha \} + (k+1) \right) f(\alpha) + \sum_{\beta \parallel \alpha} \epsilon_{\gamma\alpha} \epsilon_{\gamma\beta}  f(\beta),
	\end{align}
	where we use $\beta \parallel \alpha$ to denote that $\beta$ is a parallel neighbor of $\alpha$, and use $\epsilon_{\gamma\alpha}=\pm 1$ to denote the incidence number of $\alpha$ relative to the $(k-1)$-simplex $\gamma=\alpha \cap \beta$. More precisely, let $\alpha=\{i_0,\cdots,i_k\}$, let $\beta=\{l,i_0,\cdots,\widehat{i_j},\cdots,i_{k}\}$. $\gamma=\{i_0,\cdots,\widehat{i_j},\cdots,i_{k}\}$. We then have $\epsilon_{\gamma \alpha}=(-1)^j$ and $\epsilon_{\gamma\beta}=1$.
	
	Let $\alpha$ be a $k$-simplex of $K(\Gamma)$, define $1_{\alpha}$ to be the $k$-form such that $1_{\alpha}(\alpha)=1$ and $1_{\alpha}(\beta)=0$ for all $\alpha \neq \beta \in K_{k+1}(\Gamma)$. With a slight abuse of notation, we identify $1_{\alpha}$ with $\alpha$. $\Omega^k(\Gamma)$ can then be viewed as a $(\# K_{k+1}(\Gamma))$-dimensional real vector space with a basis $\{\alpha\}$, and $\Delta_k = \Delta|_{\Omega^k(\Gamma)}$ can be viewed as a symmetric matrix. We have
	\begin{align*}
		(\Delta_k)_{\alpha\beta}= \langle \alpha, \Delta_k (\beta) \rangle = \#\{(k+1)\textrm{-simplices}~\beta'> \alpha \} + (k+1) 
	\end{align*} 
    when $\beta = \alpha$, and
    \begin{align*}
         (\Delta_k)_{\alpha\beta}= \epsilon_{\gamma \alpha} \epsilon_{\gamma \beta}  
    \end{align*}
    when $\beta \parallel \alpha$, and $0$ for all the other cases. Let's apply the Weitzenböck decomposition to $\Delta_k$ with $R=\mathbb{R}$. We have
    \begin{align*}
   	     (\mathrm{Ric}_k)_{\alpha\beta} = \mathrm{Ric}(\Delta_k)_{\alpha\beta} = \left(\#\{(k+1)\textrm{-simplices}~\beta'> \alpha \} + (k+1) \right) \delta_{\alpha\beta} - \#\{\beta| \beta \parallel \alpha\}.
    \end{align*}    
    It is easy to see that $\mathcal{R}(\alpha)=(\mathrm{Ric}_k)_{\alpha\alpha}$ for a $k$-form $\alpha$. 
    
    To sum up, $\Delta$ can be decomposed as
    \begin{align*}
    	\Delta= B + \mathrm{Ric}
    \end{align*}
    where $B=B(\Delta)$ is the discrete counterpart of the $\nabla^*\nabla$ term in \eqref{wb}, and $\mathrm{Ric}=\mathrm{Ric}(\Delta)$ is the discrete analogue of the Weitzenböck operator.
    
	\section{Exterior covariant derivatives}\label{s_4}
    
    Let $\Gamma=(V,E)$ be a directed acyclic graph. Let $<$ be the canonical partial order on $V$ induced by the acyclic orientation on $\Gamma$.
    \begin{defn}
	Let $G$ be a Lie group. A principal $G$-bundle over $\Gamma$ is a pair $P=(V \times G, E)$. Let $(V \times V)_E = \lbrace (i,j)\in V \times V | \lbrace i,j \rbrace \in E \rbrace$. A connection on $P$ is a map $A: (V \times V)_E \rightarrow G$ such that $A(j,i) =A(i,j)^{-1}$. A gauge transformation of $P$ is a map $g: V \rightarrow G$.
	\end{defn}
	
	Let $\mathcal{A}$ be the space of connections on $P$. Let $\mathcal{G}$ be the set of gauge transformations on $P$. Note that $\mathcal{G}$ is a group and has a canonical (left) action on $\mathcal{A}$ defined by
	\begin{align*}
	(gA)(i,j) = g(i) A(i,j) g(j)^{-1}
	\end{align*}
	for $g \in \mathcal{G}$ and $A \in \mathcal{A}$. 
	\begin{rmk}\label{propgaugeaction}
		In general, the $\mathcal{G}$-action on $\mathcal{A}$ is not free. Fix a connection $A$, let 
		\[
		A_{\gamma}=A(i_1,i_2)A(i_2,i_3)\cdots A(i_{m-1},i_{m}),
		\]
		where $\gamma=i_1 \rightarrow i_2 \rightarrow \cdots \rightarrow i_m$ is a (not necessarily directed) path (or cycle) in $\Gamma$. It is easy to see that $gA=A$ if and only if
		\begin{align*}
			g(i_m)= A_{\gamma}^{-1} g(i_1) A_{\gamma}
		\end{align*} 
		for all $i_1, i_m \in V$ and all paths connecting them. 
		
		If $G$ is abelian, we have $gA=A$ if and only if $g(i)=g(j)$ for all $i, j \in V$ that can be connected by a path. In other words, the stabilizer group at each $A \in \mathcal{A}$ is isomorphic to $G$ for a connected graph. In general, if $\Gamma$ is connected, the solution $g$ to the equation $gA=A$ is uniquely determined by the value of $g$ at (any) one of the vertices of $\Gamma$. Moreover, we should have
		\begin{align*}
			g(i)=A_{lp}^{-1}g(i)A_{lp}
		\end{align*}
	    for all $i \in V$, where $lp$ is a cycle in $\Gamma$ based at $i$. It follows that the stabilizer group should be isomorphic to a subgroup of $G$ which contains the center $Z(G)$ of $G$. Note that the stabilizer group may vary for different $A \in \mathcal{A}$, and is isomorphic to $G$ for the trivial connection. 
 	\end{rmk}
    \begin{rmk}\label{propgaugeaction2}
    	 Let $\Gamma=(V,E)$ be connected. Fix a vertex $i_0 \in V$, consider the subgroup $\mathcal{G}_{i_0}$ of $\mathcal{G}$ defined by $\mathcal{G}_{i_0}=\{g \in \mathcal{G}| g(i_0)=1\}$. $\mathcal{G}_{i_0}$ acts freely on $\mathcal{A}$ since $gA=A$ implies that $g(i)=A_{\gamma}^{-1}g(i_0)A_{\gamma}=A_{\gamma}^{-1}A_{\gamma}=1$, where $\gamma$ is a path connecting $i_0$ and $i$. Note that $\mathcal{G}/\mathcal{G}_{i_0} \cong G$ acts canonically on $\mathcal{A}/\mathcal{G}_{i_0}$, making it a $G$-manifold. We have $\left(\mathcal{A}/\mathcal{G}_{i_0}\right)/G \cong \mathcal{A}/\mathcal{G}$. By Remark \ref{propgaugeaction}, the $G$-action on $\mathcal{A}/\mathcal{G}_{i_0}$ is trivial if $G$ is abelian, and we have $\mathcal{A}/\mathcal{G}_{i_0} \cong \mathcal{A}/\mathcal{G}$.
    	 
    \end{rmk}
    %\textcolor{orange}{
    %Let $\Gamma'$ be a graph equipped with a connection $A'$. Let $m: \Gamma \rightarrow \Gamma'$ be a morphism. $m$ induces a connection $A$ on $\Gamma$. 
    %}
	\begin{defn}
	Let $W$ be a vector space. We call $W_{\Gamma}=(V \times W, E)$ an associated bundle to $P$ if $W$ carries a representation of $G$ through $\rho: G \rightarrow \mathrm{GL}(W)$. A section of $W_{\Gamma}$ is a map $s: V \rightarrow W$.
	\end{defn}
	Let $s$ be a section of $W_{\Gamma}$, a gauge transformation $g$ acts on $s$ via
	\begin{align*}
	(gs)(i)= \rho(g(i))s(i).
	\end{align*}
	\begin{defn}
		A $W$-valued $k$-form on $\Gamma$ is an element in $\Omega^k(\Gamma)\otimes W$.    
	\end{defn}
    We denote the space of $W$-valued $k$-forms by $\Omega^k(\Gamma, W)$, and put $\Omega(\Gamma,W)=\bigoplus_{k=0}^{\omega(\Gamma)-1} \Omega^k(\Gamma,W)$.
	A gauge transformation $g$ should act on a $W$-valued $k$-form $f$ via
	\begin{align*}
	(gf)(i_0,\cdots,i_k) = \rho(g(i_0)) f(i_0,\cdots,i_k),
	\end{align*}
	where $i_0 < \cdots < i_k$. 
	\begin{defn}
    An $\mathrm{End}(W)$-valued $k$-form on $\Gamma$ is a map $\varphi: V \underbrace{\times \cdots \times}_{k+1} V \rightarrow \mathrm{End}(W)$ such that $\varphi(i_0,\cdots,i_k) = 0$ for all $\lbrace i_0, \cdots, i_k \rbrace \notin K_{k+1}(\Gamma)$.     
    \end{defn}
    We denote the space of $\mathrm{End}(W)$-valued $k$-forms by $\Omega^k(\Gamma, \mathrm{End}(W))$, and put $\Omega(\Gamma,\mathrm{End}(W))=\bigoplus_{k=0}^{\omega(\Gamma)-1} \Omega^k(\Gamma,\mathrm{End}(W))$. A gauge transformation $g$ should act on $\varphi$ via
    \begin{align*} 
    	(g\varphi)(i_0,\cdots,i_k)=\rho(g(i_0))\varphi(i_0,\cdots,i_k)\rho(g(i_k))^{-1}.
    \end{align*}
    \begin{rmk}
    	An $\mathrm{End}(W)$-valued $k$-form $\varphi$ is not necessarily alternating. %Therefore, $d(d\varphi)\neq 0$ in general. 
    	The alternating property is dropped because it is not preserved by gauge transformations and we want to interpret connections as $\mathrm{End}(W)$-valued $1$-forms.
    \end{rmk}
	Let $\varphi_1$ be an $\mathrm{End}(W)$-valued $p$-form and $\varphi_2$ be an $\mathrm{End}(W)$-valued $q$-form, their product $\varphi_1 \varphi_2$ is a $\mathrm{End}(W)$-valued $(p+q)$-form defined by
	\begin{align*}
		\varphi_1 \varphi_2 (i_0, \dots, i_{p+q}) = \varphi_1(i_0,\dots,i_p)\varphi_2(i_p,\dots,i_{p+q}).
	\end{align*}
    It follows that $g(\varphi_1\varphi_2)=g(\varphi_1)g(\varphi_2)$, i.e., the product is compatible with gauge transformations. We also define $[\varphi_1,\varphi_2]=\varphi_1\varphi_2 - (-1)^{pq}\varphi_2\varphi_1$.
	\begin{defn}
	Let $\varphi$ be an $\mathrm{End}(W)$-valued $p$-form and $f$ be a $W$-valued $q$-form. The action of $\varphi$ on $f$ is defined by setting
	\begin{align*}
	(\varphi f)(i_0,\cdots,i_{p+q})=\varphi(i_0,\cdots,i_p)f(i_p,\cdots,i_{p+q})
	\end{align*}
	and
	\begin{align*}
	(\varphi f)(i_{\sigma(0)},\cdots,i_{\sigma(p+q)})=(-1)^{|\sigma|}(\varphi f)(i_0,\cdots,i_{p+q})
	\end{align*}
	for $i_0 < \cdots <i_p< \cdots,i_{p+q}$, $\sigma \in S_{p+q+1}$.
	\end{defn}
	Again, one can check that $g(\varphi f)=g(\varphi)g(f)$.
	
	\begin{defn}\label{conn}
    Let $s$ be a $W$-valued $0$-form, i.e., a section of $W_{\Gamma}$. The covariant derivative $\nabla_A$ of $W_{\Gamma}$ is defined by setting
    \begin{align*}
    (\nabla_A s)(i,j)=\rho(A(i,j))s(j)-s(i),
    \end{align*}
    and 
    \begin{align*}
    (\nabla_A s)(j,i)=-(\nabla_A s)(i,j)
    \end{align*}        
    for $i<j$, $\lbrace i,j \rbrace \in E$.
    \end{defn}
    	
    %\begin{rmk}
    %Since $s(j)$ and $s(i)$ live on different vertices, to compare them, $A$ is used to ``transport" the vector $s(j)$ from $j$ to $i$, hence the name ``connection''.
    %\end{rmk}
    
    \begin{rmk}
    	We use the assumption that $\Gamma$ is a directed graph in Definition \ref{conn}. One will see later that the covariant derivative $\nabla_{A}$ can be generalized to a exterior covariant derivative $d_A$ under the assumption that $\Gamma$ is a directed acyclic graph.
    \end{rmk}
	One can check that $\nabla_{g A}(g s)=g(\nabla_{A} s)$. In fact, we have
	\begin{align*}
	\rho((gA)(i,j))(gs)(j)-(gs)(i) &= \rho(g(i))\rho(A(i,j))\rho(g(j)^{-1})\rho(g(j))s(j)-\rho(g(i))s(i) \\
	&=\rho(g(i))(\rho(A(i,j))s(j)-s(i)) \\
	&=(g\nabla_A s)(i,j).
	\end{align*}
	We also want to define $\nabla_A$ for $\mathrm{End}(W_{\Gamma})$. Let $s$ be a section of $W_{\Gamma}$. Let $\varphi$ be an $\mathrm{End}(W)$-valued $0$-form. Observe the following 
	\begin{align*}
		(\nabla_A (\varphi s) - \varphi (\nabla_A s))(i,j)&=\rho(A(i,j))\varphi(j)s(j)-\varphi(i)s(i) - \varphi(i)(\rho(A(i,j))s(j)-s(i))\\
		&=(\rho(A(i,j))\varphi(j)-\varphi(i)\rho(A(i,j)))s(j).
	\end{align*} 
    \begin{defn}
    	The covariant derivative $\nabla_A: \mathrm{End}(W_{\Gamma}) \rightarrow \Omega^1(\Gamma,\mathrm{End(W)})$ is defined by
    	\begin{align*}
    		(\nabla_A \varphi)(i,j)=\rho(A(i,j))\varphi(j)-\varphi(i)\rho(A(i,j))
    	\end{align*}
    	for $\lbrace i,j \rbrace \in E$.
    \end{defn}
    By definition, $\nabla_A$ satisfies Leibniz's rule
    \begin{align*}
    	\nabla_A(\varphi s) = \nabla_A(\varphi) s + \varphi \nabla_A(s).
    \end{align*}
    It is also easy to check that it is compatible with gauge transformations.

    \begin{defn}
    The (exterior) covariant derivative $d_A: \Omega^k(\Gamma,W) \rightarrow \Omega^{k+1}(\Gamma,W)$ is defined by
    \begin{align*}
    (d_A f)(i_0,\cdots,i_{k+1})=\rho(A(i_0,i_1))f(i_1,\cdots,i_{k+1}) + \sum_{j=1}^{k+1}(-1)^j f(i_0,\cdots,\widehat{i_j},\cdots,i_{k+1})
    \end{align*}
    and
    \begin{align*}
    (d_A f)(i_{\sigma(0)},\cdots,i_{\sigma(k+1)})=(-1)^{|\sigma|}(d_A f)(i_0,\cdots,i_{k+1})
    \end{align*}
    for $i_0 < \cdots <i_{k+1}$, where $\sigma \in S_{k+2}$.
    \end{defn}

    \begin{rmk}
        Let $d'$ denote exterior covariant derivative $d_A$ with $A=1$. It follows that
        \begin{align*}
        	(d' f)(i_0,\cdots,i_{k+1})= \sum_{j=0}^{k+1}(-1)^j f(i_0,\cdots,\widehat{i_j},\cdots,i_{k+1})
        \end{align*}
        for $i_0 < \cdots < i_{k+1}$. Note that
        \begin{align*}
        	(d'f)(i_{\sigma(0)},\cdots,i_{\sigma(k+1)})=(-1)^{|\sigma|}(d' f)(i_0,\cdots,i_{k+1})
        \end{align*}
        automatically for any $\sigma \in S_{k+2}$. Therefore, $d'$ does not depend on the partial order on $V$ and is just the exterior derivative $d$. Let $\bar{A}=\rho(A)-1 \in \Omega^1(\Gamma,\mathrm{End}(W))$. We can then write $d_A= d + \bar{A}$. 
    \end{rmk}
    
    \begin{defn}
    The (exterior) covariant derivative $d_A: \Omega^k(\Gamma,\mathrm{End}(W)) \rightarrow \Omega^{k+1}(\Gamma,\mathrm{End}(W))$ is defined by
    \begin{align*}
    (d_A \varphi)(i_0,\cdots,i_{k+1})&=\rho(A(i_0,i_1))\varphi(i_1,\cdots,i_{k+1}) + \sum_{j=1}^{k}(-1)^j \varphi(i_0,\cdots,\widehat{i_j},\cdots,i_{k+1}) \\
    &+ (-1)^{k+1}\varphi(i_0,\cdots,i_k)\rho(A(i_k,i_{k+1})).
    \end{align*}
    \end{defn}
	We have $d_A \varphi = d \varphi + [\bar{A},\varphi]$. $d_A$ reduces to $d$ for the trivial connection $A=1$. 
	\begin{prop}\label{p_4_2}
	Let $f$ be a $W$-valued $q$-form. Let $\varphi$ be a $\mathrm{End}(W)$-valued $p$-form. We have
	\begin{align*}
	d_A (\varphi f) = (d_A \varphi) f + (-1)^p \varphi (d_A f).
	\end{align*}
	\end{prop}
	
	\begin{proof}
	This follows from direct computations.
	\begin{align*}
	&(d_A \varphi) f(i_0,\cdots,i_{p+q+1}) + (-1)^p \varphi(d_A f)(i_0,\cdots,i_{p+q+1}) =  \left( \vphantom{\sum_{j=1}^{p+1}} \rho(A(i_0,i_1))\varphi(i_1,\cdots,i_{p+1}) + \right. \\
	&\left. \sum_{j=1}^{p+1}(-1)^j \varphi (i_1,\cdots,\widehat{i_j},\cdots,i_{p+1}) +  (-1)^{p+1}\varphi(i_0,\cdots,i_p)\rho(A(i_{p},i_{p+1})) \right) f(i_{p+1},\cdots,i_{p+q+1})+\\
	&(-1)^p \varphi(i_0,\cdots,i_p)\left(\rho(A(i_p,i_{p+1}))f(i_{p+1},\cdots,i_{p+q+1}) + \sum_{j=p+1}^{p+q+1} (-1)^{j-p} f (i_{p},\cdots,\widehat{i_j},\cdots,i_{p+q+1})\right) \\
	&= \rho(A(i_0,i_1))(\varphi f)(i_1,\cdots,i_{p+q+1}) + \sum_{j=1}^{p+q+1}(-1)^j (\varphi f)(i_0,\cdots,\widehat{i_j},\cdots,i_{p+q+1}) \\
	&= (d_A (\varphi f))(i_0,\cdots,i_{p+q+1})
	\end{align*}
    for $i_0 < \cdots < i_{p+q+1}$.
	\end{proof}
    
    \begin{prop}\label{p_4_1}
    	$d_A \circ d_A = F$, where $F$ is an $\mathrm{End}(W)$-valued $2$-form defined by 
    	\begin{align*}
    		F(i,j,k)=\rho(A(i,j)A(j,k))-\rho(A(i,k))
    	\end{align*}
    	We call $F$ as the curvature $2$-form of $A$. Moreover, one can show that $F=d\bar{A} + \frac{1}{2}[\bar{A},\bar{A}]=d\bar{A} + \bar{A}\bar{A}$.
    \end{prop}
    \begin{proof}
    	We have
    	\begin{align*}
    		d_A(d_Af) &= d_A(df) + d_A(\bar{A}f) 
    		= \bar{A}(df) + (d_A\bar{A})f - \bar{A}(d_A f) \\
    		&= \left(d_A\bar{A} - \bar{A}\bar{A}\right)f 
    		= \left(d \bar{A} + [\bar{A},\bar{A}] - \bar{A}\bar{A}\right)f \\
    		&= \left(d \bar{A} + \bar{A}\bar{A}\right)f,
    	\end{align*}
        where we use $d^2=0$ and the Leibniz's rule of $d_A$ to pass to the second equality. Let $F=d \bar{A} + \bar{A}\bar{A}$, we compute
        \begin{align*}
        	F(i,j,k)&=\bar{A}(j,k)-\bar{A}(i,k)+\bar{A}(i,j)+\bar{A}(i,j)\bar{A}(j,k) \\
        	&=-\bar{A}(i,k)+\rho(A(i,j))\rho(A(j,k))-1 \\
        	&=\rho(A(i,j)A(j,k))-\rho(A(i,k)).
        \end{align*}
    \end{proof}
    
    \begin{rmk}
    	%To convince the reader that 
    	Proposition \ref{p_4_1} can be also proved using direct computations. We have
    	\begin{align*}
    		&d_A(d_A f)(i_0,\cdots,i_{k+2})= \\
    		&\rho(A(i_0,i_1))(d_A f)(i_1,\cdots, i_{k+2}) + \sum_{j_1=1}^{k+2}(-1)^{j_1}(d_A f)(i_0,\cdots,\widehat{i_{j_1}},\cdots,i_{k+2}) \\
    		&=\rho(A(i_0,i_1))\left( \rho(A(i_1,i_2))f(i_2,\cdots, i_{k+2}) + \sum_{j_2=2}^{k+2} (-1)^{j_2-1} f(i_1,\cdots,\widehat{i_{j_2}},\cdots,i_{k+2}) \right) - \\
    		&(d_A f)(i_0,i_2,\cdots,i_{k+2}) +\sum_{j_1=2}^{k+2}(-1)^{j_1}(d_A f)(i_0,\cdots,\widehat{i_{j_1}},\cdots,i_{k+2}).		
    	\end{align*}
        Note that
        \begin{align*}
        	&\sum_{j_1=2}^{k+2}(-1)^{j_1}(d_Af)(i_0,\cdots,\widehat{i_{j_1}},\cdots,i_{k+2})= \sum_{j_1=2}^{k+2}(-1)^{j_1}\rho(A(i_0,i_1))f(i_1,\cdots \widehat{i_{j_1}},\cdots,i_{k+2}) + \\
        	&\sum_{1 \leq j_4 < j_1 \leq k+2}(-1)^{j_1 + j_4} f(\cdots, \widehat{i_{j_4}},\cdots,\widehat{i_{j_1}},\cdots) + \sum_{2 \leq j_1<j_4 \leq k+2}(-1)^{j_1+j_4-1}f(\cdots, \widehat{i_{j_1}},\cdots,\widehat{i_{j_4}},\cdots) \\
        	&=\sum_{j_1=2}^{k+2}(-1)^{j_1}\rho(A(i_0,i_1))f(i_1,\cdots \widehat{i_{j_1}},\cdots,i_{k+2}) + \sum_{j_1=2}^{k+2}(-1)^{j_1+1}f(i_0,\widehat{i_1},\cdots,\widehat{i_{j_1}},\cdots,i_{k+2}).
        \end{align*}
        It follows that
        \begin{align*}
    		&d_A(d_A f)(i_0,\cdots,i_{k+2})= \\
    		&\rho(A(i_0,i_1)A(i_1,i_2))f(i_2,\cdots, i_{k+2}) - \sum_{j_2=2}^{k+2} (-1)^{j_2} \rho(A(i_0,i_1) f(i_1,\cdots,\widehat{i_{j_2}},\cdots,i_{k+2}) \\
    		&-\left(\rho(A(i_0,i_2))f(i_2,\cdots,i_{k+2}) + \sum_{j_3=2}^{k+2}(-1)^{j_3-1} f(i_0,\widehat{i_1},\cdots,\widehat{i_{j_3}},\cdots,i_{k+2})\right) + \\
    		&\sum_{j_1=2}^{k+2}(-1)^{j_1}\rho(A(i_0,i_1))f(i_1,\cdots \widehat{i_{j_1}},\cdots,i_{k+2}) + \sum_{j_1=2}^{k+2}(-1)^{j_1+1}f(i_0,\widehat{i_1},\cdots,\widehat{i_{j_1}},\cdots,i_{k+2})\\
    		&=(\rho(A(i_0,i_1)A(i_1,i_2))-\rho(A(i_0,i_2)))f(i_2,\cdots,i_{k+2}) \\
    		&=F(i_0,i_1,i_2)f(i_2,\cdots,i_{k+2}) \\
    		&=(Ff)(i_0,\cdots,i_{k+2}).
    	\end{align*}
    \end{rmk}
    
    \begin{rmk}
    	Let $G$ be an abelian Lie group, for example, $\mathrm{U}(1)$. Let $P$ be a trivial principal $G$-bundle over a $n$-manifold $M$, $n \geq 2$. Let $A$ be a connection $1$-form and $F$ be its corresponding curvature $2$-form, $F = dA$, applying Stokes' theorem, we have
    	\begin{align*}
    		\int_{\gamma} A = \int_{\sigma} F,
    	\end{align*}
    	where $\sigma$ is a homology cycle of dimension $2$ on $M$ and $\gamma = \partial \sigma$. From this point of view, another reasonable definition of the curvature $2$-form on a graph is 
    	\begin{align*}
    		\tilde{F}(i,j,k)=A(i,j)A(j,k)A(k,i).
    	\end{align*}
    	 A connection is called flat if $F=0$, or equivalently, if $\tilde{F}=1$. For later use, we denote $\rho(\tilde{F})$ by $\tilde{F}_{\rho}$.
    \end{rmk}

    \begin{prop}[Second Bianchi identity]
    $d_A F=0$.
    \end{prop}
    \begin{proof}
    This follows from the observation
    \begin{align*}
    0 &= (d_A \circ d_A \circ d_A)(\cdot) - (d_A \circ d_A \circ d_A)(\cdot) \\
    &= d_A (F(\cdot)) -F (d_A(\cdot)) \\
    &= (d_A F)(\cdot),
    \end{align*}
    where we use Propositions \ref{p_4_1} and \ref{p_4_2}.
    \end{proof}
    
    \section{A generalized Weitzenböck formula}
    
    Let $\Gamma=(V,E)$ be a directed acyclic graph. Let $<$ be the canonical partial order on $V$ induced by the acyclic orientation on $\Gamma$. Let $W$ be a vector space endowed with an inner product $(\cdot,\cdot)$. Let $\rho$ be a representation of $G$ on $W$ that is orthogonal with respect to $(\cdot,\cdot)$. We put a gauge invariant inner product $\langle \cdot,\cdot \rangle$ on $\Omega(\Gamma, W_{\Gamma})$ by setting
    \begin{align*}
    	\langle f_1,f_2 \rangle = \sum_{\{i_0,\cdots,i_k\} \in K_{k+1}(\Gamma)}(f_1(i_0,\cdots,i_k),f_2(i_0,\cdots,i_k)),
    \end{align*}
    where $f_1$, $f_2$ are two $W$-valued $k$-forms. It is not hard to show that the adjoint $d^*_A$ of $d_A$ with respect to $\langle \cdot,\cdot \rangle$ takes the form
    \begin{align*}
    	d^*_A f(i_0,\cdots,i_{k-1}) = \sum_{l<i_0}\rho(A(i_0,l))f(l,i_0,\cdots,i_{k-1})+\sum_{l>i_0}f(l,i_0,\cdots,i_{k-1}).
    \end{align*}
    
    \begin{defn}
    	The connection Laplacian on $\Omega(\Gamma, W_{\Gamma})$ is defined as
    	\begin{align*}
    		\Delta_A = d_A d^*_A + d^*_A d_A.
    	\end{align*}
    \end{defn}
    
    Let $\lbrace i_0,\cdots,i_k \rbrace$ be a $(k+1)$-clique, $i_0<\cdots<i_k$. We have
    \begin{align*}
    	&d_A d^*_A f(i_0,\cdots,i_k) = \rho(A(i_0,i_1))d^*_Af(i_1,\cdots,i_k) + \sum_{j=1}^k (-1)^j  d^*_Af(i_0,\cdots,\widehat{i_j},\cdots,i_k) \\
    	&= \sum_{l<i_1}\rho(A(i_0,i_1)A(i_1,l))f(l,i_1,\cdots,i_k) + \sum_{l>i_1}\rho(A(i_0,i_1))f(l,i_1,\cdots,i_k) + \\
    	&\sum_{j=1}^k (-1)^j  \left(\sum_{l<i_0}\rho(A(i_0,l))f(l,i_0,\cdots,\widehat{i_j},\cdots,i_k) + \sum_{l>i_0}f(l,i_0,\cdots,\widehat{i_j},\cdots,i_k) \right)\\
    	&=(k+1)f(i_0,\cdots,i_k) + \sum_{l<i_1,l \neq i_0}\rho(A(i_0,i_1)A(i_1,l))f(l,i_1,\cdots,i_k) + \sum_{l>i_1}\rho(A(i_0,i_1))f(l,i_1,\cdots,i_k) \\
    	&+ \sum_{j=1}^k (-1)^j  \left(\sum_{l<i_0}\rho(A(i_0,l))f(l,i_0,\cdots,\widehat{i_j},\cdots,i_k) + \sum_{l>i_0, l \neq i_j}f(l,i_0,\cdots,\widehat{i_j},\cdots,i_k) \right).
    \end{align*}
    %where the summations inside the bracket of the last line are taken over vertices $l$ which together with $i_0,\cdots,i_{j-1}$, $i_{j+1}, \cdots,i_k$ form a $(k+1)$-cliques, $j=1,\cdots,k$, and the summations in the second last line are taken over vertices $l$ which together with $i_1, \cdots,i_k$ form a $(k+1)$-cliques.    
    We also have
    \begin{align*}
    	&d^*_A d_A f(i_0,\cdots,i_k) = \sum_{l<i_0} \rho(A(i_0,l))d_Af(l,i_0,\cdots,i_k) + \sum_{l>i_0}d_Af(l,i_0,\cdots,i_k) \\
    	&=\sum_{l<i_0<\cdots<i_k} \rho(A(i_0,l))\left(\rho(A(l,i_0))f(i_0,\cdots,i_k) - \sum_{j=0}^k(-1)^j f(l,i_0,\cdots,\widehat{i_j},\cdots,i_k)\right) + \\
    	&\sum_{i_0<l<i_1<\cdots<i_k}\left(f(i_0,\cdots,i_k)-\rho(A(i_0,l))f(l,i_1,\cdots,i_k) - \sum_{j=1}^k (-1)^j f(l,i_0,\cdots,\widehat{i_j},\cdots,i_k) \right) + \\
    	&\sum_{i_0<i_1<\cdots<l<\cdots<i_k}\left(f(i_0,\cdots,i_k)-\rho(A(i_0,i_1))f(l,i_1,\cdots,i_k)- \sum_{j=1}^k (-1)^j f(l,i_0,\cdots,\widehat{i_j},\cdots,i_k) \right) \\
    	&= \mathrm{deg}(i_0,\cdots,i_k) f(i_0,\cdots,i_k) - \sum_{l<i_0<i_1<\cdots<i_k} \rho(A(i_0,l))f(l,i_1,\cdots,i_k) - \\
    	&\sum_{i_0<l<i_1<\cdots<i_k} \rho(A(i_0,l))f(l,i_1,\cdots,i_k) - \sum_{i_0<i_1<\cdots<l<\cdots<i_k} \rho(A(i_0,i_1))f(l,i_1,\cdots,i_k) \\
    	&-\sum_{j=1}^k(-1)^j\left( \sum_{l<i_0<\cdots<i_k} \rho(A(i_0,l))f(l,i_0,\cdots,\widehat{i_j},\cdots,i_k) + \sum_{i_0<\cdots <l<\cdots<i_k}f(l,i_0,\cdots,\widehat{i_j},\cdots,i_k)\right).
    \end{align*}
    %where the summations in the last two lines are taken over vertices $l$ which together with $i_0,\cdots,i_k$ form $(k+2)$-cliques.
    
    \begin{prop}
    	For $k>0$ and $i_0<\cdots<i_k$, we have
    	\begin{align}\label{e_6_1}
    		\Delta_A f(i_0,\cdots,i_k)= \Delta^g f(i_0,\cdots,i_k) + \sum_{l<i_1} F(i_0,i_1,l)f(l,i_1,\cdots,i_k),
    	\end{align}
    	where the summation is taken over vertices $l$ which together with $i_0,\cdots,i_k$ form $(k+2)$-cliques. $\Delta^g$ can be viewed as a gauged version of the Hodge Laplacian $\Delta$, it takes the form
    	\begin{align}\label{e_6_2}
    		&\Delta^g f(i_0,\cdots,i_k) = \left(\mathrm{deg}(i_0,\cdots,i_k)+(k+1)\right) f(i_0,\cdots,i_k) \notag \\
    		&+ \sum_{l<i_1}\rho(A(i_0,i_1)A(i_1,l))f(l,i_1,\cdots,i_k) + \sum_{l>i_1}\rho(A(i_0,i_1))f(l,i_1,\cdots,i_k) \notag \\ 
    		& + \sum_{j=1}^k(-1)^j \left( \sum_{l<i_0} \rho(A(i_0,l)) f(l,i_0,\cdots,\widehat{i_j},\cdots,i_k) + \sum_{l>i_0} f(l,i_0,\cdots,\widehat{i_j},\cdots,i_k) \right),
    	\end{align}
    	where the summations in the second line are taken over vertices $l \neq i_0$ which are non-adjacent to $i_0$ and the summations inside the bracket of the third line are taken over vertices $l \neq i_j$ which are non-adjacent to $i_j$.
    \end{prop}
    
    \begin{rmk}
    	In (\ref{e_6_2}), fix $i_0$ and $i_1$, we can apply a gauge transformation $g$ such that $A(i_0,l)=1$ when $l$ is adjacent to $i_0$, and $A(i_1,l)=1$ when $l$ is adjacent to $i_1$ but non-adjacent to $i_0$. In such a gauge, we have
    	$
    		\Delta^g = \Delta
    	$
    	at the $k$-cliques containing the edge $\{i_0,i_1\}$. Note that, however, the curvature term in (\ref{e_6_1}) can not be eliminated by a gauge transformation.
    \end{rmk}

    Let's apply the Weitzenböck decomposition to $\Delta^g$. In this case, $\Delta^g$ can be viewed as a symmetric matrix whose entries are elements of the unital Banach algebra $R=\mathrm{End}(W)$. The norm $||\cdot||$ on $R$ is chosen to be the operator norm, i.e.,
    \begin{align*}
    	||r||=\sup_{||w||_W=1} ||r(w)||_W,
    \end{align*}
    where $||\cdot||_W$ is the norm on $W$ associated to the inner product $(\cdot,\cdot)$. Since we require the group action $\rho$ of $G$ to be orthogonal with respect to $(\cdot,\cdot)$, $||\rho(A(i,j))||=1$ for all $\{i,j\}\in E$, and 
    \begin{align*}
    	\mathrm{Ric}(\Delta^g)=\mathrm{Ric}(\Delta).
    \end{align*}
    \begin{prop}[Generalized Weitzenböck formula]
    	The connection Laplacian $\Delta_A$ can be decomposed as
    	\begin{align}\label{e_6_3}
    		\Delta_A = B_A + \mathrm{Ric} + F.
    	\end{align}
    	where $B_A = B(\Delta_A)$ is the discrete counterpart of the $\nabla_A^*\nabla_{A}$ term in \eqref{wbg}, the action of $F$ on a $W$-valued differential form is defined in \eqref{e_6_1}.
    \end{prop}

    \begin{rmk}
    	For $k=0$, we have
    	\begin{align}\label{connlap0}
    		\Delta_A f(i) = d^*_A d_A f(i) = \mathrm{deg}(i)f(i)-\sum_{\{l,i\} \in E} \rho(A(i,l))f(l).
    	\end{align}
    	We can apply a gauge transformation such that $\Delta_A = \Delta$ at $i \in V$. A similar formula like (\ref{connlap0}) is given in \cite{Kenyon2011} and is used there to prove a generalized version of the matrix-tree theorem in graph theory. Moreover, if we take $G$ to be $\mathbb{Z}_2$, $W$ to be $\mathbb{R}$, and the representation $\rho$ to be such that it sends $-1 \in \mathbb{Z}_2$ to $-1 \in \mathbb{R}$, we recover the notion of a signed Laplacian, which has been extensively studied in the literature. Interested readers can refer to \cite{Atay2020} for further details on this topic.
    	
    	A natural question is then if the generalized Weitzenböck formula (\ref{e_6_3}) can be applied to study problems in graph theory.
    \end{rmk}
    
    \section{Gauge invariant functionals on graphs}
     
    Let $\Gamma=(V,E)$ be a directed acyclic graph with $\omega(\Gamma) \geq 3$. Let $<$ be the canonical partial order on $V$ induced by the acyclic orientation on $\Gamma$. Let $W$ be $\mathbb{R}^n$ (or $\mathbb{C}^n$) endowed with the canonical inner product $(\cdot,\cdot)$. Let $G$ be $\mathrm{O}(n)$ (or $\mathrm{U}(n)$). Let $\rho$ be a (complex) representation of $G$ on $W$ that is orthogonal with respect to $(\cdot,\cdot)$. For an endomorphism $\varphi$, we denote its adjoint w.r.t $(\cdot,\cdot)$ as $\varphi^{\dagger}$. The gauge invariant inner product $\langle \cdot,\cdot \rangle$ on $\Omega(\Gamma,\mathrm{End}(W))$ is defined by
    \begin{align}\label{gii}
    	\langle \varphi_1,\varphi_2\rangle = \frac{1}{(k+1)!}\sum_{i_0, \cdots, i_k}  \mathrm{Tr}(\varphi_1(i_0,\cdots,i_k)\varphi_2(i_0,\cdots,i_k)^{\dagger}),
    \end{align}
    where $\varphi_1$, $\varphi_2$ are two $\mathrm{End}(W)$-valued $k$-forms. It is not hard to check that $\langle \cdot,\cdot \rangle$ is symmetric and non-degenerate, hence well-defined. Moreover, $\langle \cdot,\cdot \rangle$ is independent on the orientation of $\Gamma$ since we are summing over all the different orders of $i_0,\cdots,i_k$ in \eqref{gii}.
    
    \subsection{Yang-Mills functional}
    
    \begin{defn}
    The Yang-Mills functional on $\mathcal{A}$ is defined as
    \begin{align}\label{YM}
    \mathcal{YM}_{\rho}(A) = \frac{1}{2}\langle F, F \rangle.
    \end{align}
    \end{defn}
    \begin{prop}\label{indepord}
    	$\mathcal{YM}_{\rho}$ does not depend on the orientation of $\Gamma$.
    \end{prop}
    \begin{rmk}
    Note that 
    \begin{align*}
    F(i_0,i_1,i_2)F(i_0,i_1,i_2)^{\dagger} = 2-\tilde{F}_{\rho}(i_0,i_1,i_2)-\tilde{F}_{\rho}(i_0,i_1,i_2)^{\dagger}
    \end{align*}
    and
    \begin{align*}%\label{cycsgnf}
    \mathrm{Tr}(\tilde{F}_{\rho}(i_0,i_1,i_2))=\mathrm{Tr}(\tilde{F}_{\rho}(i_{\sigma(0)},i_{\sigma(1)},i_{\sigma(2)})^{\mathrm{sgn}(\sigma)}),
    \end{align*}
    where $\sigma \in S_3$, and
    \begin{align*}
    \tilde{F}_{\rho}(i_0,i_1,i_2)^{-1} = \tilde{F}_{\rho}(i_0,i_1,i_2)^{\dagger}.
    \end{align*}
    Therefore, we can also define the Yang-Mills functional as
    \begin{align}\label{latYM}
    \mathcal{YM}_{\rho}(A)= \frac{1}{2}\sum_{i<j<k} \mathrm{Tr}(2 - \tilde{F}_{\rho}(i,j,k) - \tilde{F}_{\rho}(i,j,k)^{\dagger}).
    \end{align}
    This form of the Yang-Mills functional can be found in a lot of lattice gauge theory literature, and is referred to as the ``Wilson action" \cite{Chatterjee2019}. Here we managed to define it on a general graph. It corresponds to the leading order term in the lattice approximation to the continuous Yang-Mills functional. See, for example, \cite{Bilson-Thompson2003} for a more detailed explanation. 
    
    The independence of \eqref{latYM} from the orientation of $\Gamma$ can be either seen from Proposition \ref{indepord}, or directly from the cyclic property of the trace.
    % and \eqref{cycsgnf}.
    \end{rmk}
    
    %\begin{rmk}
    %	\eqref{latYM} serves as an equivalent definition for the Yang-Mills functional in the discrete setting,
    %\end{rmk}
    
    Once one writes down the functional, it is very natural to ask what its Euler-Lagrange equations are. (For simplicity, we omit $\rho(\cdot)$ in relevant expressions from now on.) Note that
    \begin{align*}
    \delta F (i,j,k)=\delta A(i,j)A(j,k) - \delta A(j,k) + A(i,j) \delta A(j,k) = (d_A \delta A) (i,j,k),
    \end{align*}
    where $\delta A$ is the variation of $A$, i.e., an $\mathrm{End}(W)$-valued $1$-form.
    %, and $\delta F$ is the variation of $F$, i.e., an $\mathrm{End}(W)$-valued $2$-form. 
    We then have
    \begin{align*}
    \delta \mathcal{YM}_{\rho}(A) = \langle d_A (\delta A), F \rangle = \langle \delta A, d^*_A F \rangle,
    \end{align*}
    where $d^*_A$ is the adjoint of $d_A$ with respect to $\langle \cdot, \cdot \rangle$. Note that $\delta A$ cannot be an arbitrary $\mathrm{End}(W)$-valued $1$-form. Since $A(j,i)=A(i,j)^{-1}$, we must have 
    \begin{align*}
    \delta A(j,i)=-A(j,i)\delta A(i,j)A(j,i).
    \end{align*}
    It follows that
    \begin{align*}
    \langle \delta A, d^*_A F \rangle &= \sum_{i,j} \mathrm{Tr}(\delta A(i,j)(d^*_A F(i,j))^{\dagger}) \\
    &= \sum_{i<j} \left( \mathrm{Tr}\left(\delta A(i,j)(d^*_A F(i,j))^{\dagger}\right) + \mathrm{Tr}\left(\delta A(j,i)(d^*_A F(j,i))^{\dagger}\right) \right) \\
    &= \sum_{i<j} \mathrm{Tr}\left(\delta A(i,j)\left((d^*_A F(i,j))^{\dagger}-A(j,i)(d^*_A F(j,i))^{\dagger}A(j,i)\right)\right).
    \end{align*}
    \begin{prop}
    The Euler-Lagrange equations of the Yang-Mills functional are
    \begin{align}\label{YMeq}
    d^*_A F (i,j) = A(i,j)d^*_A F (j,i)A(i,j)
    \end{align}
    for all $\lbrace i,j \rbrace \in E$. We also refer to \eqref{YMeq} as the Yang-Mills equations.
    \end{prop}
    
    It is not hard to work out an explicit formula for $d^*_A$, which is
	\begin{align*}
	d^*_A \varphi(i_0,\cdots,i_{k-1})&=\frac{1}{k+1} \sum_{l} \left( \vphantom{\sum_{j=1}^{k-1}} A(i_0,l)\varphi(l,i_0,\cdots,i_{k-1}) \right.\\
	&\left. + \sum_{j=1}^{k-1} (-1)^j \varphi(i_0,\cdots,i_{j-1},l,i_j,\cdots,i_{k-1}) + (-1)^k \varphi(i_0,\cdots,i_{k-1},l)A(l,i_{k-1})\right).
	\end{align*}
	With this formula in hands, it is not hard to show that
	\begin{align*}
	d^*_A F (i_0,i_1) = -\sum_l F(i_0,l,i_1).
	\end{align*}
	Thus, the Yang-Mills equations (\ref{YMeq}) become
	\begin{align*}
	\sum_{l} F(i,l,j)=A(i,j)\sum_{l} F(j,l,i)A(i,j).
	\end{align*}
	Or equivalently,
	\begin{align*}
	\sum_l \tilde{F}_{\rho}(i,l,j)= \sum_l \tilde{F}_{\rho}(i,j,l).
	\end{align*}

    Note that the Yang-Mills is bounded from both below and above. More precisely, we have $0 \leq \mathcal{YM}_{\rho}(A) \leq 2 n \omega_3$, where $\omega_3$ is the number of $3$-cliques in $\Gamma$, and $n$ is the (complex) dimension of $W$.
    \begin{prop}
    	%For any (connected) graph $\Gamma$ with $\omega(\Gamma) \geq 3$, t
    	The global minimum and the global maximum of the Yang-Mills functional can always be achieved.
    \end{prop}
    \begin{proof}
        The global minimum can be achieved by the trivial connection, or more generally, by a flat connection $A$. The global maximum can be achieved by setting $A(i,j)=-1$ for all $\{i,j\} \in E$, or more generally, by a connection $A$ satisfying $\tilde{F}_{\rho}=-1$.
    \end{proof}
    \begin{defn}
    	We call a solution $A$ to the Yang-Mills equations trivial if $A= \pm 1$ up to a gauge transformation.
    \end{defn}
	
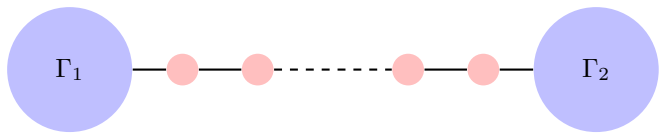
\begin{figure}[!h]
    	\centering
    	\begin{tikzpicture}[shorten >=1pt,->]
    		\tikzstyle{vertex}=[circle,fill=red!25,minimum size=12pt,inner sep=2pt]
    		\tikzstyle{vertexg}=[circle,fill=blue!25,minimum size=12pt,inner sep=12pt]
    		\node[vertexg] (G_1) at (0,0) {$\Gamma_1$};
    		\node[vertex] (G_2) at (1.5,0) {};
    		\node[vertex] (G_3) at (2.5,0) {};
    		\node[vertex] (G_4) at (4.5,0) {};
    		\node[vertex] (G_5) at (5.5,0) {};
    		\node[vertexg] (G_6) at (7,0) {$\Gamma_2$};
    		\draw[thick] (G_1) -- (G_2)-- cycle;
    		\draw[thick] (G_2) -- (G_3) -- cycle;
    		\draw[thick,dashed] (G_3) -- (G_4) -- cycle;
    		\draw[thick] (G_4) -- (G_5) -- cycle;
    		\draw[thick] (G_5) -- (G_6) -- cycle;
    	\end{tikzpicture}
    	\caption{A graph $\Gamma$ obtained by connecting two graphs $\Gamma_1$ and $\Gamma_2$ with a path.}
    	\label{decomp}
    \end{figure}

    Since the functional (\ref{YM}) is gauge invariant, we can choose a proper gauge to simplify the calculations. Let $\Gamma$ be a connected graph. Let $T=(V_t,E_t)$ be a spanning tree of $\Gamma$, i.e., a subgraph that is a tree which includes all of the vertices of $\Gamma$. One can then consider the following procedure. 
    \begin{enumerate}
    \item Pick a vertex $a \in V_t$, let $g(a)=1$.
    \item For any vertex $b$ adjacent to $a$, choose $g(b)=A(a,b)$.
    \item For any vertex $c$ adjacent to $b$, choose $g(c)=A(c,b)$. Since $T$ is a tree, $c$ is distinct from $a$ and its adjacent vertices $b$, $g(c)$ is well defined.
    \item Repeat the Step 3 until all vertices of $T$ are exhausted. 
    \end{enumerate}
    We then have $A(i,j)=1$ for all $\lbrace i,j \rbrace \in E_t$. We refer to this procedure as the spanning tree gauge fixing. Its generalization to disconnected graphs is called as the spanning forest gauge fixing.
    \begin{rmk}
    By applying the spanning forest gauge fixing to $\Gamma$, one can easily see that $\mathcal{A}/\mathcal{G}=\mathrm{pt}$ if $\Gamma$ is a forest. In other words, gauge theory is only sensitive to non-acyclic graphs, i.e., graphs with cycles.
    \end{rmk}
    \begin{rmk}
    	One may wonder how to find a global maximizer of the Yang-Mills functional in the spanning tree gauge. Let $T$ be a spanning tree of $\Gamma$. For any two vertices $i$ and $j$ of $T$, we set $A(i,j)=(-1)^{d_{ij}-1}$ if $\{i,j\}$ is an edge of $\Gamma$, where $d_{ij}$ is the length of the path connecting $i$ and $j$ in $T$. In particular, $A(i,j)=1$ if $\{i,j\}$ is an edge of $T$, i.e., we are in the spanning tree gauge of $\Gamma$. For any such triangle $\{i,j,k\}$ in $\Gamma$, we assume that, without loss of generality, $d_{ik}=d_{ij}+d_{jk}$. It follows that
    	\begin{align*}
    		\tilde{F}_{\rho}(i,j,k)=A(i,j)A(j,k)A(k,i)=(-1)^{d_{ij}-1}(-1)^{d_{jk}-1}(-1)^{-d_{ik}+1}=-1.
    	\end{align*}
    	Therefore, $A$ is a global maximizer of the Yang-Mills functional.
    \end{rmk}
	\begin{exmp}
	Let $G=\mathrm{U}(1)$.
	\begin{figure}[!h]
	\centering
	\begin{tikzpicture}[shorten >=1pt,->]
    \tikzstyle{vertex}=[circle,fill=red!25,minimum size=12pt,inner sep=2pt]
    \node[vertex] (G_1) at (0,0) {1};
    \node[vertex] (G_2) at (1,1.732)   {2};
    \node[vertex] (G_3) at (2,0)  {3};
    \draw[color=blue,thick] (G_2) -- (G_1) -- (G_3) -- cycle;
    \draw[thick] (G_2) -- (G_3) -- cycle;
    \end{tikzpicture}
    \caption{The complete graph $K_3$ with a spanning tree.}
    \label{K3}
    \end{figure}
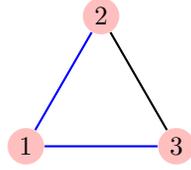
	We choose the spanning tree $T$ as in Figure \ref{K3}. There is only one degree of freedom, namely, $A(2,3)$. The solutions to (\ref{YMeq}) are the trivial solutions
	\begin{align*}
	A(2,3)= \pm 1.
	\end{align*}
	\end{exmp}
	%In fact, the global minimum and maximum of the Yang-Mills functional can always be achieved. Using the spanning tree gauge fixing, it is not hard to show that both the minimizer and the maximizer are unique up to gauge transformations. We refer to the corresponding solutions as the trivial solutions.
    \begin{exmp}
	Let $G=\mathrm{U}(1)$.
	\begin{figure}[!h]
	\centering
	\begin{tikzpicture}[shorten >=1pt,->]
    \tikzstyle{vertex}=[circle,fill=red!25,minimum size=12pt,inner sep=2pt]
    \node[vertex] (G_1) at (0,0) {1};
    \node[vertex] (G_2) at (1,1.732)   {2};
    \node[vertex] (G_3) at (2,0)  {3};
    \node[vertex] (G_4) at (1,0.57735)  {4};
    \draw[color=blue,thick] (G_1) -- (G_2) -- cycle;
    \draw[color=blue,thick] (G_1) -- (G_3) -- cycle;
    \draw[color=blue,thick] (G_1) -- (G_4) -- cycle;
    \draw[thick] (G_2) -- (G_4) -- cycle;
    \draw[thick] (G_2) -- (G_3) -- (G_4) -- cycle;
    \end{tikzpicture}
    \caption{The complete graph $K_4$ with a spanning tree.}
    \label{K4}
    \end{figure}
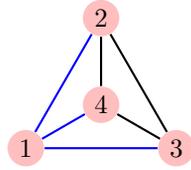
	We choose the spanning tree $T$ as in Figure \ref{K4}. There are three degrees of freedom, namely, $A(2,3)$, $A(2,4)$, $A(3,4)$. Let $G=\mathrm{U}(1)$. The nontrivial solutions to (\ref{YMeq}) are
	\begin{align*}
	&A(2,3)=A(2,4)= A(3,4)= \pm 1, \\
	&A(2,3)=A(2,4)^{-1}=-A(3,4)^{-1}= \exp(i\alpha), \\
	&A(2,3)=-A(2,4)=A(3,4)= \exp(i\alpha), \\
	&A(2,3)=-A(2,4)=-A(3,4)^{-1}= \exp(i\alpha).
	\end{align*}  
    \begin{figure}[!h]
    	\centering
    	\begin{tikzpicture}
    		\draw[thick] (2,2) circle (1.5cm);
    		\draw[thick] (2,2) ellipse (1.5cm and 1cm);
    		\draw[thick] (2,2) ellipse (1.5cm and 0.5cm);
    	\end{tikzpicture}
    	\caption{The space of nontrivial solutions to the $\mathrm{U}(1)$ Yang-Mills equations on $K_4$.}
    	\label{solK4}
    \end{figure}
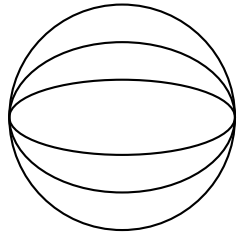
    The space of nontrivial solutions is depicted in Figure \ref{solK4}. The two intersection points correspond to the solutions
    \begin{align*}
    	A(2,3)=-A(2,4)=A(3,4)=\pm i.
    \end{align*}
    This space has a $S_3$-symmetry, which is inherited from the $S_3$-symmetry of the pair $(K_4,T)$.
	\end{exmp}
    \begin{rmk}
    	The computation of solutions to the $\mathrm{U}(1)$ Yang-Mills equations on the complete graph $K_n$ becomes much more harder for $n \geq 5$. However, we know the space of nontrivial solutions (if nonempty) should have a $S_{n-1}$-symmetry.
    \end{rmk}
    To simplify the computation, a natural idea is to decompose the graph into easily computable pieces. For example, if a (connected) graph $\Gamma$ can be obtained by connecting two graphs $\Gamma_1$ and $\Gamma_2$ connected by a path (see Figure \ref{decomp}), then it is easy to show that the space of solutions to the Yang-Mills equations on $\Gamma$ is the Cartesian product of the spaces of solutions to the Yang-Mills equations on $\Gamma_1$ and $\Gamma_2$.

    \subsection{Yang-Mills-Higgs functional}
    
    \begin{defn}
    	Let $\Gamma$ be a graph with $\omega(\Gamma) \geq 3$. The Yang-Mills-Higgs functional is defined as
    	\begin{align}\label{YMH}
    		\mathcal{YMH}_{\rho}(A,\phi) = \frac{1}{2}\langle F, F \rangle + \frac{1}{2}\langle d_A \phi, d_A \phi \rangle + V(\phi),
    	\end{align}
        where $\phi$ is a section of $W_{\Gamma}$ and $V$ is an non-negative function on the space of sections of $W_{\Gamma}$.
    \end{defn}
    Let's derive the Euler-Lagrange equations to \eqref{YMH}. For simplicity, we set $V(\phi)=0$. The variation of the second term of \eqref{YMH} with respect to $A$ is 
	\begin{align*}
		\langle \delta A \phi, d_A \phi \rangle &= \sum_{i<j} \mathrm{Tr}( \delta A \phi(i,j)(d_A\phi(i,j))^{\dagger}) \\
		&= \sum_{i<j} \mathrm{Tr}\left(\delta A(i,j)\phi(j)(A(i,j)\phi(j)-\phi(i))^{\dagger}\right) \\
		&= \sum_{i<j} \mathrm{Tr}\left(\delta A(i,j)\left(\phi(j)\phi(j)^{\dagger}A(j,i)-\phi(j)\phi(i)^{\dagger}\right)\right).
	\end{align*}
    Recall that
    \begin{align*}
    	\langle \delta A, d^*_A F \rangle = \sum_{i<j} \mathrm{Tr}\left(\delta A(i,j)\left((d^*_A F(i,j))^{\dagger}-A(j,i)(d^*_A F(j,i))^{\dagger}A(j,i)\right)\right).
    \end{align*}
    We arrive at the Euler-Lagrange equations for the connection $A$, which are
    \begin{align}
    	d^*_A F(i,j)-A(i,j)(d^*_A F(j,i))A(i,j)=\phi(i)\phi(j)^{\dagger}-A(i,j)\phi(j)\phi(j)^{\dagger},
    \end{align}
    or equivalently,
    \begin{align}
    	\sum_l \tilde{F}_{\rho}(i,l,j)-\sum_l \tilde{F}_{\rho}(i,j,l) = A(i,j)\phi(j)\phi(j)^{\dagger}A(j,i) - \phi(i)\phi(j)^{\dagger}A(j,i).
    \end{align}
    The variation of the second term of \eqref{YMH} with respect to $\phi$ is 
    \begin{align*}
    	\langle d_A \delta \phi, d_A \phi \rangle = \langle \delta \phi, \Delta_A \phi(i) \rangle = 0,
    \end{align*}
    where $\Delta_A$ is the connection Laplacian. By \eqref{connlap0}, the Euler-Lagrange equations for the scalar field $\phi$ are
    \begin{align}
    	\mathrm{deg}(i)\phi(i)-\sum_l A(i,l)\phi(l)=0.
    \end{align}

	\section{Conclusions and future directions}
	
	In this paper, we have developed a discrete setting for gauge theories which preserves most of the flavor of the original continuous setting. Possible future research directions can be a more detailed study of the solutions to the Yang-Mills(-Higgs) equations, and incorporations of further concepts and theorems (e.g., Wilson's area law \cite{Wilson1974}) from lattice gauge theory into this framework.
    
    \section*{Acknowledgement}
    
    The author would like to thank Jürgen Jost for many helpful discussions. This work was supported by the International Max Planck Research School Mathematics in the Sciences.

\end{document}